\documentclass[11pt]{amsart}
\usepackage{amsmath,amssymb,amscd}
\usepackage{eucal}

\newtheorem{thm}{Theorem}[section]
\newtheorem{prop}{Proposition}[section]
\newtheorem{lemma}{Lemma}[section]
\newtheorem{cor}{Corollary}[section]
\newtheorem{conj}{Conjecture}[section]

\theoremstyle{remark}
\newtheorem{remark}{Remark}[section]

\theoremstyle{definition}
\newtheorem{definition}{Definition}[section]

\newcommand{\nn}{\nonumber}
\newcommand{\bea}{\begin{eqnarray*}}
\newcommand{\eea}{\end{eqnarray*}}
\newcommand{\bean}{\begin{eqnarray}}
\newcommand{\eean}{\end{eqnarray}}
\newcommand{\ben}{\begin{equation}}
\newcommand{\een}{\end{equation}}
\newcommand{\be}{\begin{equation*}}
\newcommand{\ee}{\end{equation*}}
\newcommand{\half}{\frac{1}{2}}

\def\ol{\overline}

\def\ora{\overset{\rightarrow}}

\newcommand{\bF}{{\mathbb F}}     
\newcommand{\bN}{{\mathbb N}}     
\newcommand{\bP}{{\mathbb P}}     
\newcommand{\bQ}{{\mathbb Q}}     
\newcommand{\bZ}{{\mathbb Z}}     

\newcommand{\cL}{{\mathcal L}}

\newcommand{\cR}{{\mathcal R}}

\newcommand{\cP}{{\mathcal P}}
\newcommand{\cU}{{\mathcal U}}

\newcommand{\cW}{{\mathcal W}}

\newcommand{\nop}{\ \Xi\ }

\begin{document}

\title[A proof of GV conjecture for local toric Calabi-Yau]
{A simple proof of Gopakumar--Vafa conjecture for local toric
Calabi-Yau manifolds}

\author{Pan Peng}
\dedicatory{Department of Mathematics, UCLA}

\address{Department of Mathematics, UCLA}
\email{ppeng@math.ucla.edu}

\date{10/18/2004}

\begin{abstract}
We prove Gopakumar-Vafa conjecture for local toric Calabi-Yau
manifolds. It is also proved that the local Goparkumar-Vafa
invariants of a given class at large genus vanish.
\end{abstract}

\maketitle


\section{Introduction}

Let $X$ be a Calabi--Yau threefold. $K_{\Sigma }^{g} (X)$ is
defined as the genus $g$, 0 marked point Gromov--Witten invariant
of $X$ in the curve class $\Sigma \in H_{2} (X,\bZ )$. From
M--theory, Gopakumar and Vafa \cite{GV2} express the invariants
$K_{\Sigma }^{g} (X)$ in terms of integer invariants $N_{\Sigma
}^{g} (X)$ obtained by BPS state counts . In the physics
literature, these invariants should be the Euler characteristics
of moduli space of D-brane. Unfortunately, this moduli space has
not been defined so far in mathematics. We will take the following
Gopakumar--Vafa formula as the definition of Gopakumar-Vafa
invariant $N_\Sigma^g(X)$.
\ben%
 \sum_{\Sigma \neq 0}\sum_{g\geq 0}
    K_\Sigma^g(X) t^{2g-2}e^{-\omega \cdot \Sigma}
 = \sum_{\Sigma \neq 0}\sum_{n\geq 1} \sum_{g\geq 0}
    \frac{N_\Sigma^g(X)}{n}
    (2\sin(\frac{nt}{2}))^{2g-2}e^{-n\omega\cdot \Sigma}
    \label{eqn: GV-invariants definition formula}
\een%
One can recursively obtain Gopakumar--Vafa invariants from
Gromov-Witten invariants through the above formula and vice versa.
It's not hard to see that Gopakumar-Vafa invariants are rational
number from the expression. The following Gopakumar--Vafa
conjecture captures our intensive interests.
\begin{conj}{\bf(Gopakumar--Vafa)}
$\forall \Sigma \in H_2(X,\bZ)$ and $g\geq 0$, $N_\Sigma^g\in
\bZ$.
\end{conj}
From (\ref{eqn: GV-invariants definition formula}), there is no
reason (at least mathematically) that these $N_\Sigma^g$'s should
be integers.

In general, in order to obtain Gopakumar-Vafa invariants from
Gromov-Witten invariants, one has to get a way to compute the
corresponding topological string amplitudes. For lower genus, one
may use mirror symmetry or localization method. However, if one
goes to higher genus, the calculation become more and more
complicated. Neither of these two methods is very practical. This
situation has changed since a beautiful discovery of the duality
between Chern--Simons theory and Gromov-Witten theory \cite{GV3}.
Chern-Simons theory provides a very efficient way to compute all
genus amplitudes for toric Calabi-Yau manifolds.

Recently, in \cite{LLLZ} J. Li, C.-C. Liu, K. Liu and J. Zhou
build the mathematical foundation of topological vertex theory,
which allows us to rigorously write down all genus amplitudes for
any toric Calabi-Yau manifold directly from its toric diagram.
Thus it may provide us a new way to approach Gopakumar--Vafa
conjecture. For physics literature about topological vertex, one
may refer \cite{AKMV}.

In this paper, we consider the local toric Calabi-Yau geometry,
which is the geometry of Fano surfaces within Calabi-Yau
threefold. The corresponding all genus A-model amplitudes obtained
from its relation to Chern-Simons amplitudes can be organized in a
simple way. In local Calabi-Yau, a tubular neighborhood of the
Fano surface is equivalent to the total space of the canonical
bundle. Note that local Calabi-Yau manifolds are noncompact.

Let $S$ be a nonsingular, projective, toric, Fano surface and
denote by $K_S$ its canonical bundle. The Gromov-Witten theory of
local Calabi-Yau geometry of $S$ is defined by excess integral.
Let $\overline{M}_{g}(S, \Sigma)$ denote the moduli space of
stable maps from connected genus g curves with 0 marked point to
$S$ representing the class $\Sigma\in H_2(S, \bZ)$. The local
$(\Sigma, g)$ Gromov-Witten invariant of $S$ is defined by
\bea%
 K_\Sigma^g(S) = \int_{[\overline{M}_{g}(S, \Sigma)]^{vir}} e(R^1 \pi_* ev^* K_S)
\eea%
Where $\pi: \cU\rightarrow \ol{M}_g(S, \Sigma)$ is the universal
curve, $ev: \ol{M}_g(S, \Sigma)\rightarrow S$ is the evaluation
map. We still take the formula (\ref{eqn: GV-invariants definition
formula}) as the definition of local Gopakumar--Vafa invariant
$N_\Sigma^g(S)$ while we use the corresponding local
Gromov--Witten invariant $K_\Sigma^g(S)$ instead.

However, the invariant $\widehat{N}^g_\Sigma(S)$ we mention later
as local Gopakumar--Vafa invariant defers slightly from the above.
The relation between $\widehat{N}^g_\Sigma(S)$ and $N^g_\Sigma(S)$
is $\widehat{N}^g_\Sigma(S)=(-1)^{g-1}N_\Sigma^g(S)$. We call them
both local Gopakumar-Vafa invariants.

In this paper, we proved that for any local toric Calabi-Yau
manifold, the Gopakumar-Vafa conjecture is true.  More precisely,
given any $\Sigma$ as above, we have:

{\bf Theorem 5.2} \ \ \  $P^{(S)}_\Sigma(x)=\sum_{g\geq 0}
\widehat{N}^g_\Sigma(S) x^g\in \bZ[x]$.

Before beginning the proof we may note that, from our result,
Gopakumar--Vafa invariants vanish at large genera. This is an very
interesting phenomenon which is pointed out in \cite{GV2}. After
our calculations on Hirzebruch surfaces and $\bP^2$, we find that
the degree of $P_\Sigma$ is actually the arithmetic genus
$A_\Sigma$ of the curve with the degree $\Sigma$. Moreover,
$\widehat{N}^{A_\Sigma}_\Sigma=\chi(L(\Sigma))$, where $L(\Sigma)$
is the holomorphic line bundle associated to the devisor $\Sigma$.
All these seems to imply that Gopakumar--Vafa invariants should be
related to the embedding curves.

Our proof of Theorem 5.2 contains two parts. First, we extract
Gopakumar--Vafa invariants by proving a generalized M\"obius
inversion formula ({\bf Theorem \ref{thm:mobius}}). Thus we obtain
an explicit formula for $P_\Sigma(x)$, while r.h.s of the formula
merely contains the Chern-Simons invariants of Hopf links. The
standard calculation of Schur function then leads to the following
property: $P_\Sigma(x)=\frac{1}{c}h(x)$, where $c\in\bN$, $h(x)\in
\cL[x]$ ($\cL[x]$ is defined by (\ref{eqn:definition of L[x]})).
Here comes a crucial observation: $P_\Sigma(x)$ is a polynomial.
In the second part, we show a pattern theorem ({\bf Theorem
\ref{thm:main-tech}}). Theorem \ref{thm:main-tech} asserts if
$P_\Sigma(x)$ satisfies some conditions which exactly fits our
situation, then $P_\Sigma(x)\in \bZ[x]$. We prove Theorem
\ref{thm:main-tech} by deliberately pairing the data when multiple
covering happens. This pairing shows $c=1$, which is equivalent to
$P_\Sigma\in \cL[x]$. $P_\Sigma(x)\in\bZ[x]$ then follows from the
structure of $\cL[x]$.

In our proof, the pattern theorem (Theorem \ref{thm:main-tech})
doesn't depend on the concrete expression of all genus amplitudes
of local toric Calabi-Yau manifolds. The technique we develop here
may extend to general toric Calabi-Yau using the topological
vertex theory. Now this work is in its preparation.

It is worth mentioning that a new idea of understanding
Gopakumar-Vafa conjecture \cite{LiLZ} comes out recently. They
regard topological string amplitudes as equivariant indices and
translate Gopakumar-Vafa conjecture into a formula of an infinite
product, where the integrality naturally follows.

The rest of this paper is organized as follows. We fix some
notations on partition in section 2. Then we prove some results on
symmetric functions and Chern-Simons invariants of the Hopf link
in section 3. In section 4, we give an explicit formula for
Gopakumar-Vafa invariants. Our main result is proved in section 5.
Some examples are discussed in section 6.

{\bf Acknowledgements}. {\em I would like to thank my advisor,
Professor Kefeng Liu, for his  constant encouragement and a lot of
inspiring discussions. I am very grateful to Xiaowei Wang for his
many helpful comments and friendship. I also would like to thank
Professor Jian Zhou for pointing out some misleading parts and
some mistakes. At last, I want to thank Professor Amer Iqbal for
his explaining me the physical reason of the vanishing property of
Gopakumar-Vafa invariants at large genera.}

\section{Preliminaries}

\subsection{Partitions}
A partition is a sequence of non-negative integers
$$\lambda=(\lambda_1, \lambda_2, \cdots, \lambda_r, \cdots)$$
such that
$$\lambda_1\geq\lambda_2\geq\cdots\geq\lambda_r\geq\cdots$$
and containing only finitely many non-zero terms. The degree of
$\lambda$ is defined by
$$|\lambda|=\lambda_1+\lambda_2+\cdots$$
If $|\lambda|=d$, we denote it by $\lambda\vdash d$. The length of
$\lambda$ is defined by $l(\lambda):=Card\{j:\lambda_j\neq 0 \}$.
The number $m_i(\lambda)=Card\{j: \lambda_j=i \}$ is called the
multiplicity of $i$ in $\lambda$. We may rewrite a partition
$\lambda$ as $(1^{m_1(\lambda)}2^{m_2(\lambda)}\cdots
n^{m_n(\lambda)}\cdots)$. One can identify a partition with its
Young diagram or $\{ (i, j)\in \bZ_+^2: j\leq\lambda_i \}$. We
will freely refer them when needed as long as without any
confusion. By transposing the Young diagram of $\lambda$, one gets
the corresponding partition $\lambda^t$. $\forall x=(i,j) \in
\lambda$, let
\begin{align*}
c(x)=j-i, \ \  h(x)=\lambda_i+\lambda_j^t-i-j+1, \ \
n(\lambda)=\sum_i (i-1)\lambda_i.
\end{align*}

For a partition $\lambda$, it's easy to see that
$|Aut\lambda|=\prod_{i\geq 1}m_i(\lambda)!$. Define the following
constants associated to the partition $\lambda$:
\begin{align*}
u_\lambda
    & = \frac{l(\lambda)!}{\prod_i m_i(\lambda)!}  =
    \frac{l(\lambda)!}{|Aut\lambda|}\\
z_\lambda
    &= |Aut\lambda| \cdot \prod_{i\geq 1}\lambda_i \\
k_\lambda
    &=  \sum_i \lambda_i (\lambda_i -2i +1) \\
\theta_\lambda
    &= \frac{(-1)^{l(\lambda)-1}(l(\lambda)-1)!}{|Aut\lambda|}
\end{align*}
and let $u_{(0)}=\theta_{(0)}=z_{(0)}=1$. Note that $k_\lambda = -
k_{\lambda^t}$ and $k_\lambda$ is even.

\begin{lemma}
$\sum_{x\in\lambda} h(x) - \frac{k_\lambda}{2} = 2n(\lambda) +
|\lambda|$.
\end{lemma}

\begin{proof}
\bea%
 k_\lambda
    &=& \sum_i\lambda_i (\lambda_i -2i +1) \\
    &=& 2\sum_i\frac{\lambda_i (\lambda_i -1)}{2} - 2\sum_i \lambda_i
        (i-1) \\
    &=& 2 (n(\lambda^t)-n(\lambda))
\eea%
Follow (\cite{M} page 11), one has %
\bea%
 \sum_{x\in\lambda}h(x) &=& n(\lambda)+n(\lambda^t)+|\lambda|
\eea%
Then direct calculation leads to the proof.
\end{proof}

\subsection{Infinite series and partition}

Given a sequence $x=(x_1, x_2, \cdots, x_n, \cdots)$, let
$|x|=\sum_{i\geq 1}x_i$. For a partition $\lambda$, define
$x_\lambda=\prod_{i\geq 1}x_{\lambda_i}$, $x_{(0)}\equiv 1$.

\begin{lemma}\label{lemma:f(|x|)}
Let $f(t)=\sum_{n\geq 0} a_n t^n$, $x$ is as above. Then
$$
f(|x|) = \sum_{|\lambda|\geq 0} a_{l(\lambda)} x_\lambda u_\lambda
$$
\end{lemma}

\begin{proof}
Note that $|x|^d=\sum_{l(\lambda)=d}x_\lambda u_\lambda$. The
proof then follows from a simple calculation.
\end{proof}

Apply the above lemma to $\exp(x)$ and $\log(x)$, we have
\begin{cor}
\begin{eqnarray}
\exp \left( {\sum_{n\geq 1}b_n x^n}\right) & = & 1+
\sum_{|\lambda|\geq
1}\frac{b_\lambda x^{|\lambda|}}{|Aut\lambda|} \label{eqn:exp}\\
\log \left(1+\sum_{n\geq 1}b_n x^n \right) & = & \sum_{|\lambda|
\geq 1}b_\lambda x^\lambda \theta_\lambda \label{eqn:log}
\end{eqnarray}
\end{cor}

\section{Symmetric functions}

\subsection{Basis}
Let $x=(x_1, x_2, \cdots )$ be a sequence, $\lambda$ a partition.
Denote by $Per(\lambda)$ = \{ $\delta$: $\delta$ is a
permutation of $\lambda$ \}. Follow the traditional notations:%
\begin{align*}%
m_\lambda
    &= \sum_{\delta \in Per(\lambda)} x^\delta, &
p_n
    &= m_{(n)}, &
h_n
    &= \sum_{\lambda\vdash n}m_\lambda%
\end{align*}%
Let $p_\lambda=\prod_i p_{\lambda_i}$. Similarly define
 $h_\lambda$. Schur function $S_\lambda=\det(h_{\lambda_i -i+j})_{1\leq i, j \leq
 l(\lambda)}$.
 It's easy to see that%
\begin{align*}
H(x,t) &= \sum_{r\geq 0} h_r t^r = \prod_{i} (1-x_i t)^{-1}
\end{align*}

\subsection{Chern-Simons invariants of the Hopf link $\cW_{\lambda\mu}(q)$}

Define:
$$ \cW_\lambda(q)=q^{\frac{k_\lambda}{4}}
\prod_{1\leq i< j\leq d}\frac{[\lambda_i-\lambda_j+j-i]}{[j-i]}
\prod_{i=1}^{l(\lambda)} \prod_{v=1}^{\lambda_i}
\frac{1}{[v-i+l(\lambda)]}$$ and
$$\cW_{\lambda\mu}(q)=\cW_\lambda(q) S_\mu
(q^{\lambda_1-\half}, q^{\lambda_2-\frac{3}{2}}, \cdots,
q^{\lambda_n-n+\half}, \cdots).$$ $\cW_\lambda(q)$ can be
simplified as
$$\cW_\lambda(q)=q^{\frac{k_\lambda}{4}}\prod_{x\in \lambda} \frac{1}{[h(x)]}
=S_\lambda(q^{-\half}, q^{-\frac{3}{2}}, \cdots, q^{-n+\half},
\cdots)
$$
where $[m]=q^{\frac{m}{2}}-q^{-\frac{m}{2}}$. Hence%
\ben \label{eqn:wR1R2}
 \cW_{\lambda\mu}(q)=q^{-\frac{\lambda+\mu}{2}}
 S_\lambda(q^\alpha) S_\mu(q^{\lambda+\alpha})
\een%
where $q^\alpha=(1, q^{-1}, q^{-2}, \cdots)$. It is proved in \cite{Z2} that%
\ben \label{eqn:sym-w-lambda-mu}
 \cW_{\lambda^t \mu^t}(q)= (-1)^{|\lambda|+|\mu|} \cW_{\lambda \mu}(\frac{1}{q})
\een%
Follow (\cite{M} page 27)
\bea%
 H(q^{\alpha},t)
 &=& \prod_{i=1}^\infty (1-q^{-i+1}t)^{-1}  \\
 &=& \sum_{r=0}^\infty \frac{t^r}{\phi_r(\frac{1}{q})} \\
 &=& \sum_{r=0}^\infty
 \frac{\phi_r(\frac{1}{q})q^{\frac{r(r+1)}{2}}}{\prod_{i=1}^r
 [i]^2}t^r
\eea%
where $\phi_r(q)=\prod_{i=1}^r (1-q^i)$. Then correspondingly
\bean%
 h_r(q^\alpha)=\frac{\phi_r(\frac{1}{q})
 q^{\frac{r(r+1)}{2}}}{\prod_{i=1}^r [i]^2}\label{eqn:h-r(q-a,t)}
\eean%
Using the results of $H(q^{\alpha},t)$, we can calculate
$H(q^{\lambda+\alpha}, t)$ in the following way:
\bea%
 H(q^{\lambda+\alpha}, t)
 &=& \prod_{i=1}^\infty (1-q^{\lambda_i-i+1}t)^{-1} \\
 &=& \prod_{i=1}^{l(\lambda)} \frac{1-q^{-i+1}t}{1-q^{\lambda_i
    -i+1}t} H(q^\alpha, t) \\
 &=& \prod_{i=1}^{l(\lambda)} \frac{1-q^{-i+1}t}{1-q^{\lambda_i -i+1}t}
    \sum_{r=0}^\infty \frac{\phi_r(\frac{1}{q})q^{\frac{r(r+1)}{2}}}{\prod_{i=1}^r [i]^2}
\eea%
Observe that one can expand $(1-q^{\lambda_i-i+1}t)^{-1}$ and get
the coefficient of $t^j$, which is an element of $\bZ[q,q^{-1}]$.
After collecting the $t^r$ terms in the above formula, we have
\ben%
 h_r(q^{\lambda+\alpha}) = \frac{\rho(q)}{\prod_{i=1}^r [i]^2}\label{eqn:h-r(qlambda+alpha)}
\een%
Here $\rho(q) \in \bZ[q, q^{-1}]$. Recall
$S_\lambda=\det(h_{\lambda_i-i+j})_{1\leq i, j \leq l(\lambda)}$.
From Lemma \ref{lemma:f(q)of-beta}, we know $[n]^2$ is a
polynomial of $[1]^2$ with integer coefficients. Combine
(\ref{eqn:h-r(q-a,t)}) and (\ref{eqn:h-r(qlambda+alpha)}), we
obtain:

\begin{lemma}\label{lemma:the form of Schur function of q-alpha}
$S_\lambda(q^\alpha), \ S_\mu(q^{\lambda+\alpha})$ are of form
$\frac{a(q)}{b([1]^2)}$, where $a(q)\in \bZ[q,q^{-1}]$, $b(x)\in
\bZ[x]$ and is monic.
\end{lemma}

Let $\cP$ denote the space of all partitions. Define%
\bea
 \cW_{\overrightarrow{R}, \gamma}=\cW_{R^NR^1} \cW_{R^1R^2}\cW_{R^2R^3}\cdots
 \cW_{R^{N-1}R^N}  (-1)^{\sum_{i=1}^N \gamma_i |R^i|}
 q^{\half \sum_{i=1}^N k_{R^i}\gamma_i}
\eea%
where $\overrightarrow{R}=(R^1, R^2, \cdots, R^N) \in \cP^N$ and
$\gamma=(\gamma_1, \gamma_2, \cdots, \gamma_N)\in \bZ^N$. Denote
by $\| \ora{R} \| = \sum_{i=1}^N |R^i|$ and $\ora{R^t}= ((R^1)^t,
(R^2)^t, \cdots, (R^N)^t )$.
Let $R^{N+1}=R^1$, we rewrite $\cW_{\overrightarrow{R}, \gamma}$ as%
\bean \label{eqn:w-R-gamma}%
 \cW_{\overrightarrow{R}, \gamma} &=&
 q^{-\sum_{i=1}^N |R^i|}
    \prod_{i=1}^N S_{R^i}(q^\alpha)
    \prod_{i=2}^{N+1}S_{R^i}(q^{R^{i-1}+\alpha})\\
  && \times (-1)^{\sum_{i=1}^N \gamma_i |R^i|}
    q^{\half \sum_{i=1}^N k_{R^i}\gamma_i} \nn
\eean%
Notice that $k_\lambda$ is even for any partition $\lambda$, from
(\ref{eqn:w-R-gamma}) and lemma \ref{lemma:the form of Schur
function of q-alpha}, we have

\begin{cor}\label{cor:W-R-Gamma}
$\cW_{\overrightarrow{R}, \gamma}$ is of form
$\frac{a(q)}{b([1]^2)}$, where $a(q)\in \bZ[q, q^{-1}]$, $b(x)\in
\bZ[x]$ and is monic.
\end{cor}

\begin{lemma} \label{lemma:W-Rt-q}
$\cW_{\ora{R^t}, \gamma}(q)=\cW_{\ora{R}, \gamma}(\frac{1}{q})$.
\end{lemma}

\begin{proof}
Before proving the lemma, let's look at the following statement
first.

\begin{lemma}\label{lemma:f(q)of-beta}
$f(q)\in \bZ[q, q^{-1}]$, $f(q)=f(q^{-1})$. Then $f\in \bZ[x]$,
where $x=[1]^2$.
\end{lemma}

\begin{proof}
$f(q)=f(q^{-1})$, and $f\in \bZ[q, q^{-1}]$. $f$ has to be of the
form $\sum_{n=0}^N a_n (q^n+q^{-n})$ for some $N$. It's sufficient
to show $\forall n$, $q^n+q^{-n}\in \bZ[x]$.

$n=1$, $q^1+q^{-1}=x + 2$. $n=k$. If $k$ odd,
$q^n+q^{-n}=(q+q^{-1})(q^{k-1}+q^{k-2}+ \cdots + q^{-(n-1)})$; if
$k$ even, $q^n+q^{-n}=(q^{\frac{n}{2}}+q^{\frac{-n}{2}})^2-2$. So
by induction, the proof is completed.
\end{proof}

From lemma \ref{lemma:f(q)of-beta} and corollary
\ref{cor:W-R-Gamma}, the proof of lemma \ref{lemma:W-Rt-q} follows
easily from the formula (\ref{eqn:sym-w-lambda-mu}).
\end{proof}

\section{Local Gopakumar-Vafa invariants}
\subsection{A generalized M\"obius inversion formula}

Let $\mu(n)$ be the M\"obius function. That is:
\be%
\mu(n)=
\left\{%
\begin{array}{ll}
    (-1)^r, & n=p_1\cdot p_2 \cdots p_r \ \
    \hbox{and} \ \ p_i\neq p_j,\ \
    if \ \  i\neq j ; \\
    0, & \hbox{otherwise.} \\
\end{array}%
\right.
\ee%

\begin{lemma}\label{lemma:mobius funtion delta-1n}
$\sum_{k|n}\mu(k)=\delta_{1n}$
\end{lemma}

\begin{proof}
$n=1$, trivial. If $n>1$, say $n=\prod_{i=1}^d p_i^{r_i}$, $p_i$'s
are distinct and
$r_i>0$. Then we have%
\bea
 \sum_{k|d} \mu(k)=\sum_{i=0}^d \begin{pmatrix}
    d \\ i \end{pmatrix}
    (-1)^i = (1-1)^d = 0.
\eea%
\end{proof}

Let $\Sigma=(k_1, k_2, \cdots, k_n)\in \bZ^n_{\geq 0}$.  Denote by
$n|\Sigma$ if $n|k_i, \ \forall 1\leq i \leq n$. For any $r\in
\bQ$, let  $r\Sigma \triangleq(r\cdot k_1, r \cdot k_2, \cdots,
r\cdot k_n) \in \bQ^n$.

\begin{thm}\label{thm:mobius}
Let $f$, $g$ be two functions defined on $\bZ_{\geq 0}^n \times
\bZ_{\geq 0}$. $\alpha$ is a function on $\bZ_+$ such that
$\alpha(n)\alpha(m)=\alpha(nm)$. If
\ben%
 \sum_{n|\Sigma}\alpha(n)f(\frac{\Sigma}{n},nm)=g(\Sigma,m)
 \label{eqn:formula of mobius invertion in the lemma}
\een%
holds for any $\Sigma \in \bZ_{\geq 0}^n$, $m \in \bZ_{\geq 0}$.
Then we have
\ben%
 f(\Sigma,m)=\sum_{n|\Sigma}\mu(n)\alpha(n)g(\frac{\Sigma}{n},nm).
 \label{eqn:solution of mobius inversion}
\een%
\end{thm}

\begin{proof}
If $\alpha \equiv 0$, it's trivial. Otherwise,
$\alpha(n)\alpha(m)=\alpha(nm)$ implies $\alpha(1)=1$. Note that
given $g(\Sigma,m)$'s, the solution of $f(\Sigma, m)$ in
(\ref{eqn:formula of mobius invertion in the lemma}) is unique.
Hence it's sufficient to show (\ref{eqn:solution of mobius
inversion}) is a solution of (\ref{eqn:formula of mobius invertion
in the lemma}).
\bea%
 \sum_{n|\Sigma}\alpha(n)f(\frac{\Sigma }{n},nm)
&=& \sum_{n|\Sigma}\alpha(n)\sum_{k|\frac{\Sigma}{n}}\mu(k)a(k)
    g(\frac{\Sigma}{kn},knm) \\
&=& \sum_{kn|d}\alpha(n)\alpha(k)\mu(k)g(\frac{\Sigma}{kn},knm) \\
&=& g(\Sigma,m) + \sum_{p>1}\sum_{kn=p}\alpha(nk)\mu(k)
    g(\frac{\Sigma}{kn},knm) \\
&=& g(\Sigma,m) +\sum_{p>1}\alpha(p)g(\frac{\Sigma}{p},pm)
    \sum_{k|p}\mu(k) \\
 &=& g(\Sigma,m)
\eea%
Where in the last step, we've used Lemma \ref{lemma:mobius funtion
delta-1n}. The proof is completed.
\end{proof}

\subsection{Explicit formula for local Gopakumar-Vafa invariants}
 Given a nonsingular, projective, toric, Fano Surface $S$, the topological
string amplitudes for the local toric Calabi-Yau manifold $K_S$
(\cite{Z2}\cite{AMV}\cite{I}\cite{EK}) is: \bea
 Z_{top\ str}^{(S)}
 &=& \sum_{R^1, R^2, \cdots, R^N} \cW_{R^NR^1}
    \cW_{R^1R^2} \cW_{R^2 R^3} \cdots \cW_{R^{N-1}R^N}
    (-1)^{\sum_{i=1}^N\gamma_i |R^i|} \\
 && \times q^{\half \sum_{i=1}^N k_{R^i} \gamma_i} \exp\{ -\sum_{i=1}^N t_i
    |R^i|\} \\
 &=& \sum_{
   \overrightarrow{R}=(R^1, R^2, \cdots, R^N)} \cW_{\ora{R}, \gamma}
   \exp \left( -\sum_{i=1}^N t_i |R^i| \right)
\eea%
Where $\gamma_i$ is the self-intersection numbers of the rational
curves associated to the i-th edge, $t_i$'s are linear
combinations of K\"ahler parameters. In our calculation, $\gamma$
is fixed. We may omit $\gamma$ and write $\cW_{\ora{R}}$ instead
of $\cW_{\ora{R}, \gamma}$.

\bea
 F_{closed}
 &=& \sum_{\Sigma\in H_2(S, \bZ),\ \Sigma \neq 0} \ \ \sum_{n=1}^\infty
    \sum_{g=0}^\infty \frac{\widehat{N}^g_\Sigma [n]^{2g-2}}{n}
    e^{-n \omega\cdot \Sigma}\\
 &=& \sum_{\Sigma\in H_2(S,\bZ),\ \Sigma\neq 0} e^{-\omega\cdot \Sigma} \sum_{n| \Sigma}
    \frac{1}{n} \sum_{g=0}^\infty \widehat{N}^g_{\frac{\Sigma}{n}} [n]^{2g-2}
\eea%
Where $\widehat{N}^g_\Sigma$'s are the local Gopakumar-Vafa
invariants.

Fix a basis in $H_2(X, \bZ)$, one can identify $H_2(X,
\bZ)=\bZ^m$. We only consider $\Sigma=(k_1, k_2, \cdots, k_m) \in
\bZ^m_{\geq 0}$. Denote by $|\Sigma|=\sum_{i=1}^m k_i$,
$\tilde{P}_{\Sigma}(x)= \sum_{g=0}^\infty \widehat{N}^g_\Sigma
x^{g-1}$. Let $\omega=(\omega_1, \omega_2, \cdots, \omega_m)^t$.
Here we write $\omega$ as an $m\times 1$ matrix, $\Sigma$ a
$1\times m$ matrix. The relation between $T\triangleq (t_1, t_2,
\cdots, t_N)^t$ and $\omega$ is $T=A\cdot \omega$, where $A$ is an
$N\times m$ matrix. Denote by $|\ora{R}|=(|R_1|, |R_2|, \cdots,
|R_N|)$.

Then%
\bea
 F_{closed}
 &=& \log \left( 1+\sum_{\|\ora{R}\| \geq 1} \cW_{\ora{R}}
    e^{-|\ora{R}|\cdot T}
    \right) \\
 &=& \log \left( 1+\sum_{d=1}^\infty \eta_d \right) \\
 &=& \sum_{\lambda} \theta_\lambda
 \eta_\lambda
\eea%
Where $\eta_d= \sum_{|\Sigma|=d}\sum_{|\ora{R}|\cdot A=\Sigma}
\cW_{\ora{R}} \exp\left( - |\ora{R}|\cdot T \right)$,
$\eta_\lambda=\prod_{i\geq 1} \eta_i$.

Combine the above two formulas,%
\bea
 \sum_{n|\Sigma} \frac{1}{n} \tilde{P}_{\frac{\Sigma}{n}}([n]^2)
 &=& \sum_{\lambda\vdash |\Sigma|} \theta_\lambda
 \langle\eta_{\lambda}(q)\rangle_{\Sigma}
\eea%
Where $\langle \eta_\lambda \rangle_\Sigma$ is the coefficient of
$e^{-\omega \cdot \Sigma}$ in $\eta_\lambda$. Let $P_\Sigma(x)=x
\tilde{P}_\Sigma(x)=\sum_{g\geq 0}\widehat{N}^g_\Sigma x^g$.
By Theorem \ref{thm:mobius}, from the above formula we can derive%
\bean
 P_\Sigma([1]^2)=[1]^2 \cdot \sum_{n|\Sigma}\frac{\mu(n)}{n}
    \sum_{\lambda\vdash \frac{|\Sigma|}{n}} \theta_\lambda
    \langle\eta_{\lambda}(q^n)\rangle_{\frac{\Sigma}{n}} \label{eqn:explicit-formula}
\eean%

We fix the following notation in this paper
\ben%
 \cL[x]=\left\{
\frac{a(x)}{b(x)}: a(x),\ b(x) \in \bZ[x], b(x)\neq 0 \ \hbox{and}
\ b(x)\  \hbox{is monic} \right\} \label{eqn:definition of L[x]}
\een%

\begin{prop}\label{prop: L[x] is a ring}
$\cL[x]$ is a ring with conventional addition and multiplication.
\end{prop}

\begin{proof}
Trivial calculation.
\end{proof}

By Corollary \ref{cor:W-R-Gamma}, Lemma \ref{lemma:W-Rt-q} and
Proposition \ref{prop: L[x] is a ring}, we have:

\begin{lemma} $d\in \bN$,
$\langle\eta_d(q)\rangle_{\Sigma} \in \cL[x]$.
\end{lemma}

\begin{prop}\label{prop: eta-lambda in the ring L[x]}
$\lambda$ is a partition, let $l=l(\lambda)$, then
$$\langle\eta_{\lambda}(q)\rangle_{\Sigma} = \sum_{\Sigma_1+\Sigma_2+\cdots +
\Sigma_{l}=\Sigma}\langle\eta_{\lambda_1}(q)\rangle_{\Sigma_1}\cdots
\langle\eta_{\lambda_{l}}(q)\rangle_{\Sigma_{l}}$$ Then by the
above lemma, $\langle\eta_{\lambda}(q)\rangle_{\Sigma} \in \cL[x]$
\end{prop}

\begin{proof}
By direction calculation.
\end{proof}

\section{The integrality of local Gopakumar-Vafa invariants}

\subsection{Pattern theorem}

$r\in\bQ$, $p$ is a prime number. Let $r=p^k \frac{m}{n}$, where
$m,\ n,\ k\in \bZ$, $p\nmid mn$. Define $\xi_p(r)=k$. Define
$$\cR[x]=\left\{\frac{c(x)}{d(x)}:\ c(x),\ d(x)\in\bZ[x]
\right\}$$

\begin{definition}
$\forall a(x),\ b(x)\in \cR[x]$. Let
$$a(x)-b(x)= r\frac{f(x)}{g(x)}$$ where $r\in\bQ$, $f$
and $g$ are primitive polynomials in $\bZ[x]$. Let $p$ be a prime
number, $k\in\bZ$. Denote by
\bea%
 a(x)\equiv b(x) \mod(p^k)
\eea%
if  $\xi_p(r)\geq k$; denote by
\bea%
 a(x)\nop b(x) \mod(\frac{1}{p})
\eea%
if $\xi_p(r)\geq 0$.
\end{definition}

\begin{remark}\label{remark:mod 1/p}
If
$
 a(x)\nop b(x)\mod(\frac{1}{p}),\ \
 b(x)\nop c(x)\mod(\frac{1}{p})
$. It's easy to see that $a(x)\nop c(x)\mod(\frac{1}{p}).$
\end{remark}

\begin{lemma}\label{devident}
$a$, $b$ $\in \bZ$, $r\in \bN$, $p$ is a prime number, $f(x)\in
\bZ[x]$. The following statements hold:
\begin{enumerate}
    \item[(A)] $\left.\frac{a}{\gcd(a,b)}\right|
        \begin{pmatrix} a\\ b \end{pmatrix}$.
    \item[(B)] $p^r|(a^{p^r}-a^{p^{r-1}})$.
    \item[(C)] $f(x)\in \bZ[x]$, $(f(x))^{p^r}\equiv (f(x^p))^{p^{r-1}}\mod (p^r)$
    \item[(D)] $g(x)\in\cL[x]$, $g(x)^{p^r}-g(x^p)^{p^{r-1}}=p^r\cdot
        h(x)$, where $h(x)\in \cL[x]$
\end{enumerate}
\end{lemma}

\begin{proof}
(A). $\begin{pmatrix} a \\ b \end{pmatrix} = \frac{a}{b}
\begin{pmatrix}
  a-1 \\
  b-1
\end{pmatrix}
\Rightarrow \left. \frac{a}{\gcd (a,b)}\right|
\begin{pmatrix}
  a \\
  b \\
\end{pmatrix}$.

(B). If $a$ is $p$, the claim is obviously true since $p^{r-1}>r$.
If $\gcd(a,p)=1$, by Fermat theorem, $a^{p-1}\equiv 1 \mod(p)$.
Then
\bea%
 a^{p^r}-a^{p^{r-1}}
 &=& a^{p^{r-1}}((kp+1)^{p^{r-1}(p-1)}-1)\\
 &=& a^{p^{r-1}}\left(\sum_{i=1}^{p^{r-1}}(kp)^i\left( \begin{array}{c}
   p^{r-1} \\
   i \\
 \end{array} \right) \right)\\
 &\equiv& 0 \mod (p^r)
\eea In the last step, we used (A) and the fact $i>\log(i)$.

(C). By (A) and (B), $(f(x))^{p^r}-(f(x^p))^{p^{r-1}}\equiv
(f(0))^{p^r} - (f(0))^{p^{r-1}} \mod(p^r)\equiv 0 \mod (p^r)$.

(D). $g(x)=\frac{a(x)}{b(x)}$, where $b(x)$ is monic. %
\bea
 g(x)^{p^r}-g(x^p)^{p^{r-1}}
 &=& \frac{a(x)^{p^r}b(x^p)^{p^{r-1}}-a(x^p)^{p^{r-1}}
    b(x)^{p^r}}{b(x)^{p^r}b(x^p)^{p^{r-1}}}\\
\eea%
By (C):%
\bea%
 && a(x)^{p^r}b(x^p)^{p^{r-1}}-a(x^p)^{p^{r-1}}b(x)^{p^r} \\
 & \equiv & a(x^p)^{p^{r-1}}b(x^p)^{p^{r-1}}
    -a(x^p)^{p^{r-1}}b(x)^{p^r} \mod(p^r) \\
 & \equiv & a(x^p)^{p^r}b(x)^{p^{r}}
    -a(x^p)^{p^{r-1}}b(x)^{p^r} \mod(p^r) \\
 & \equiv & 0 \mod(p^r)
\eea%
(D) is then proved.
\end{proof}

With these preparation, we have the following pattern theorem:

\begin{thm}\label{thm:main-tech}
$f_\Sigma (x), \ g_\Sigma^{(d)}(x), \ \psi_n(x)$ satisfy the
following conditions:
\begin{enumerate}
    \item $f_\Sigma(x)\in\bQ[x]$.
    \item $g_\Sigma^{(d)}(x)\in \cL[x], \ \forall d$. Let
        $\lambda=(\lambda_1, \cdots, \lambda_{l(\lambda)})$ be a
        partition, define
        $$g_\Sigma^{(\lambda)}=\sum_{\Sigma_1+\Sigma_2+\cdots
            +\Sigma_{l(\lambda)}=\Sigma}
            g_{\Sigma_1}^{(\lambda_1)}(x)\cdots
            g_{\Sigma_{l(\lambda)}}^{(\lambda_{l(\lambda)})}(x)
        $$
    \item $\psi_k(x)$'s are monic polynomials of degree $k$ such
    that $\psi_{np}(x)-\psi_n(x)^p\equiv 0 \mod(p)$, where $p$ is prime.
    \item The following equation holds:
        $$
         f_\Sigma(x)=x\cdot \sum_{n|\Sigma}\frac{\mu(n)}{n}
         \sum_{\lambda\vdash \frac{|\Sigma|}{n}} \theta_\lambda
         g_{\frac{\Sigma}{n}}^{(\lambda)}(\psi_n(x))
        $$
\end{enumerate}
Then, $f_\Sigma(x)\in \bZ[x]$.
\end{thm}

\begin{proof}
Denote by $H_\Sigma(x)=\sum_{\lambda\vdash |\Sigma|}
\theta_\lambda g_{\Sigma}^{(\lambda)}(x)$. Before proceeding the
proof, let's look at the following lemma.

\begin{lemma}\label{lemma:tech}
Notations are as above, then for any prime number $p$, we have
$$ H_{p\Sigma}(x)-\frac{1}{p}H_\Sigma (\psi_p(x)) \nop
 0\mod \left(\frac{1}{p} \right).$$
\end{lemma}

\begin{proof}
By (A) in Lemma \ref{devident}, $\theta_\lambda$ is not an integer
if and only if there exists $k>1$ such that $k|m_i(\lambda), \
\forall i$. However, by condition (2), the terms in $H_\Sigma$ are
not in $\cL[x]$ will be when $\lambda=(\lambda_1^{n_1}\cdots
\lambda_k^{n_k})$ and the corresponding term has to be $\pm
\frac{m}{l} (g_{\Sigma_1}^{(\lambda_1)})^{n_1} \cdots
(g_{\Sigma_k}^{(\lambda_k)})^{n_k}$, where $m$ is some integer,
$l=\gcd(n_1,\cdots, n_k)$. Otherwise $l$ will be cancelled by the
same term with permutations in $H_\Sigma$. Obviously, for this
term, $n_1\Sigma_1+\cdots +n_k \Sigma_k=\Sigma$. Then $l|\Sigma$.
So $H_\Sigma=h(x)+ \frac{a(x)}{c\cdot b(x)}$, where $h(x),\
\frac{a(x)}{b(x)} \in\cL[x]$, $c|\Sigma$.

First, let's consider $p>2$.  Then $p$ is odd. We only consider
$1/p$ factor in the terms of $H_{p\Sigma}$. The number of terms in
$H_{p\Sigma}$ which contains factor $1/p$ is equal to that of
$H_\Sigma$. We pair them in a way that the general term will be:
\bean
 A &\triangleq&
    (-1)^{p^{r-1}\cdot a \cdot s}
    \frac{1}{p^r\cdot a}\left\{ \zeta\cdot
    \left( (g_{\Sigma_1}^{(\lambda_1)}(x))^{s_1}
    \cdots (g_{\Sigma_k}^{(\lambda_k)}(x))^{s_k}\right)^{p^r\cdot a}
    \label{eqn:A-pr}\right.
    \\
 && \left. - \eta\cdot\left( (g_{\Sigma_1}^{(\lambda_1)}(\psi_p(x)))^{s_1}
    \cdots (g_{\Sigma_k}^{(\lambda_k)}(\psi_p(x)))^{s_k}
    \right)^{p^{r-1}\cdot a}
    \nn
    \right\}
    \\
 &\nop& (-1)^{p^{r-1}\cdot a \cdot s}
    \frac{\eta}{p^r\cdot a}
    \left( \phi(x)^{p^r}-\phi(x^p)^{p^{r-1}}\right)\mod(\frac{1}{p})
    \label{eqn:A-second-step}
\eean
 where $a,\  b\in \bZ$ and $\gcd(p,a)=1$ and
 \bea%
  \zeta &=&
        \begin{pmatrix}
            p^ras \\
            p^ras_1, p^ras_2, \cdots, p^ras_k
        \end{pmatrix}\\
  \eta &=& \begin{pmatrix}
            p^{r-1}as \\
            p^{r-1}as_1, p^{r-1}as_2, \cdots, p^{r-1}as_k
        \end{pmatrix}\\
  \phi(x) &=&
        (g_{\Sigma_1}^{(\lambda_1)}(x))^{as_1}
    \cdots (g_{\Sigma_k}^{(\lambda_k)}(x))^{as_k}
\eea%
 In
 (\ref{eqn:A-second-step}), we have used Theorem \ref{thm:number theory}, $\psi_p(x)\equiv x^p \mod(p)$ and
 $p^r|p^{p^{r-1}}$.
 By (D) in lemma
 \ref{devident}, we have $A \nop 0 \mod(\frac{1}{p})$.

If $p=2$ and $2|\Sigma$, similar to (\ref{eqn:A-pr}), we still get
the subtraction. Then we can still use (4) in lemma \ref{devident}
to prove $A \nop 0 \mod(\frac{1}{2})$.

If $p=2$ and $2\nmid \Sigma$. Different from (\ref{eqn:A-pr}), it
will be the summation of the two terms.  But we only have to show
$A \nop 0 \mod(\half)$ with at most $\frac{1}{2}$ appearing, which
will be obtained by $2|(1+1)$.

The lemma is proved.
\end{proof}

Now, we come back to the proof of the theorem. By the condition
(4),%
\bea
 f_\Sigma(x)=x\cdot \sum_{n|\Sigma}\frac{\mu(n)}{n}
         \sum_{\lambda\vdash \frac{|\Sigma|}{n}} \theta_\lambda
         g_{\frac{\Sigma}{n}}^{(\lambda)}(\psi_n(x)) \label{eqn:f-sigma-x}
\eea%
We have $f_\Sigma(x)=\frac{a(x)}{c\cdot b(x)}$, where $a(x), \
b(x) \in \bZ[x]$, $c|\Sigma$ and $b(x)$ is monic. To see
$c|\Sigma$, we only notice that in the above formula, $n|\Sigma$
and in the term $\sum_{\lambda\vdash \frac{|\Sigma|}{n}}
\theta_\lambda g_{\frac{\Sigma}{n}}^{(\lambda)}(\psi_n(x))$, if
there is any term which is of form $\frac{1}{d}h(x)$, where
$h(x)\in \cL[x]$, follow the first paragraph of the proof of Lemma
\ref{lemma:tech}, $d|\frac{\Sigma}{n}$.

We will show $f_\Sigma(x)\in \cL[x]$ by proving: $\forall
p|\Sigma$ prime, $f_\Sigma(x)\nop 0 \mod(\frac{1}{p})$.

\bea%
 f_\Sigma(x)
 &=& x\cdot \sum_{n|\Sigma}\frac{\mu(n)}{n}
    \sum_{\lambda\vdash \frac{|\Sigma|}{n}} \theta_\lambda
    g_{\frac{\Sigma}{n}}^{(\lambda)}(\psi_n(x)) \\
 &=& x \sum_{n|\Sigma,\  p\nmid n} \frac{\mu(n)}{n}
    H_{\frac{\Sigma}{n}}(\psi_n(x))
    + x \sum_{n|\Sigma,\  p\nmid n} \frac{\mu(pn)}{pn}
    H_{\frac{\Sigma}{pn}}(\psi_{pn}(x)) \\
 &\nop& x\cdot \sum_{n|\Sigma, \ p\nmid n}\frac{\mu(n)}{n}
    \left( H_{\frac{\Sigma}{n}}(\psi_n(x))-\frac{1}{p}
    H_{\frac{\Sigma}{pn}}((\psi_{n}(x))^p) \right) \mod(\frac{1}{p}) \\
 & \nop & 0 \mod (\frac{1}{p})
\eea%
Where in the third step, we used the condition (3) and
$\mu(pn)=-\mu(n)$; in the last step, lemma \ref{lemma:tech} is
applied.

Then $f_\Sigma(x)=\frac{a(x)}{b(x)}\in \cL[x]$. $f_\Sigma(x)\in
\bQ[x]$, so we can write $f_\Sigma(x)=\frac{1}{n}h(x)$, where $n
\in \bN$, $h(x)\in \bZ[x]$. Assume $n\nmid h(x)$, then $n>1$. On
the one hand,
$$
 a(x)=\frac{1}{n}b(x)h(x).
$$
Since $n\nmid h(x)$, we have $n|b(x)$. On the other hand, $b(x)$
is monic, so $b(x)$ is primitive, which is a contradiction. Then
$n|h(x)$ and $f_\Sigma(x)\in \bZ[x]$. The proof is completed.
\end{proof}

\subsection{Main results}

\begin{thm} $P_{\Sigma}(x)=\sum_{g=0}^\infty
\widehat{N}^g_{\Sigma} x^g \in \bZ[x]$.
\end{thm}

\begin{proof}
By (\ref{eqn:explicit-formula})%
\bea%
 P_\Sigma([1]^2)=[1]^2\sum_{n|\Sigma}\frac{\mu(n)}{n}
    \sum_{\lambda\vdash \frac{|\Sigma|}{n}} \theta_\lambda
    \langle\eta_{\lambda}(q^n)\rangle_{\frac{\Sigma}{n}}
\eea%
Let $x=[1]^2=(q^{\half}-q^{-\half})^2$, $\psi_n(x)=[n]^2$. First,
we show $P_\Sigma(x)$ satisfies (1) of the Theorem
\ref{thm:main-tech}. Note that $P_\Sigma(x)=\frac{a(x)}{b(x)}$,
where $a(x), b(x)\in\bZ[x]$. In fact,
$b(x)=q^k\prod_{i_s}(1-q^{i_s})^2$ is a finite product for some
$k$ and $i_s$'s, which can be obtained from (\ref{eqn:h-r(q-a,t)})
and (\ref{eqn:h-r(qlambda+alpha)}). As a Laurent series in $q$ for
$0<|q|<1$, $q=0$ is a pole of $\frac{a(x)}{b(x)}$. If there are
infinitely many $g$'s such that $\widehat{N}^g_\Sigma$ is nonzero
for the fixed $\Sigma$, then $q=0$ is an essential singular point
of $P_\Sigma(x)$, which is a contradiction. So $P_\Sigma(x)$ is a
polynomial. $\widehat{N}^g_\Sigma$'s are rational numbers, which
can be easily obtained from the definition. $P_\Sigma(x)\in
\bQ[x]$.

Then by Proposition \ref{prop: eta-lambda in the ring L[x]}, one
finds that our situation is exactly like the Theorem
\ref{thm:main-tech} as long as we show the functions $\psi_n$'s
have the same properties as (3) in the Theorem
\ref{thm:main-tech}.

By the definition of $\psi_n$, we only
need to show $\psi_p(x)-x^p\equiv 0\mod (p)$.%
\bea%
 [1]^{2p} &=& \sum_{k=1}^{p-1} \left( \begin{array}{c}
   2p \\
   k \\
 \end{array}
    \right) (-1)^k (q^{p-k}+q^{-(p-k)})+
    (q^p+q^{-p}-2)+ \left( \begin{array}{c}
      2p \\
      p \\
    \end{array} \right)+2 \\
 &=& \sum_{k=1}^{p-1} \left( \begin{array}{c}
   2p \\
   k \\
 \end{array}
    \right) (-1)^k (q^{n-k}+q^{-(n-k)}) + [p]^2 + \left( \begin{array}{c}
      2p \\
      p \\
    \end{array} \right)+2
\eea%
equivalently, %
\bea%
 [p]^2= -\sum_{k=1}^{p-1} \left( \begin{array}{c}
   2p \\
   k \\
 \end{array}
    \right) (-1)^k (q^{n-k}+q^{-(n-k)}) + [1]^{2p} - \left( \begin{array}{c}
      2p \\
      p \\
    \end{array} \right)-2
\eea%
 By lemma \ref{lemma:f(q)of-beta}, $q^n+q^{-n}\in \bZ[x]$, so
$[p]^2\in \bZ[x]$. By (A) of the Lemma \ref{devident},
$p|\begin{pmatrix} 2p \\ k
\end{pmatrix}$ when $1\leq k \leq p-1$.  $[1]^{2p}\rightarrow 0$ as
$[1]^2\rightarrow 0$ implies as a polynomial of $x$, $[p]^2$ has 0
as the constant term. Then $\psi_p(x)-x^p=[p]^2 - [1]^{2p} \equiv
0 \mod(p)$.

The proof is completed.
\end{proof}

\begin{remark}
Actually, if one wants, one can precisely calculate the degree of
$P_\Sigma(x)$ using formula (\ref{eqn:explicit-formula}) and the
following formula about Schur function.
$$
 S_\lambda(x)=\sum_{|\rho|=|\lambda|}
    \frac{\chi_\lambda(\rho)}{z_\rho}p_\rho(x).
$$
where $\chi_\rho$ is the character of representation $\rho$.
\end{remark}

\section{Some examples}

We may leave the results of Hirzebruch surfaces $\bF_k$ and
$\bP^2$ as the following theorem. We use the same notation as
those in the introduction: $P_{\Sigma}^{(S)}(x)=\sum_{g\geq 0}
\widehat{N}^g_\Sigma(S) x^g$, where $\widehat{N}^g_\Sigma(S)$ is
the local Gopakumar-Vafa invariant of the given class $\Sigma\in
H_2(S, \bZ)$ for the Fano surface $S$. Let $A_\Sigma(S)$ be the
arithmetic genus of the curve with the degree $\Sigma\in H_2(S,
\bZ)$. $L(\Sigma)$ denotes the holomorphic line bundle associated
to the devisor $\Sigma$.

\begin{thm}
We have
\begin{enumerate}
 \item $\deg(P_d^{(\bP^2)}(x))=A_d(\bP^2) =\half (d-1)(d-2)$. \\
    $\deg(P_{(d,m)}^{(\bF_k)}(x))=A_{(d,m)}(\bF_k) =(d-1)(m-1-\frac{kd}{2})$, \\
 \item Moreover,
 \bea%
  \widehat{N}_d^{A_{d}(\bP^2)}(\bP^2) &=&(-1)^d\chi(L(\Sigma))\\
  &=& \frac{(-1)^d}{2}
  (d+1)(d+2),\\
  \widehat{N}_{(d,m)}^{A_{(d,m)}(\bF_k)}(\bF_k)
  &=&(-1)^{dk}\chi(L(\Sigma))\\
  &=& (-1)^{dk}(d+1)(m+1-\frac{kd}{2}).
\eea%
\end{enumerate}
\end{thm}

\begin{proof}
We do computation for $\bP^2$ as an example. Let $H\in H_2(\bP^2,
\bZ)$ be the generator, it satisfies $H\cdot H=1$. So the
topological string amplitudes will be
\bea%
 Z_{top\ str}^{(\bP^2)}
  &=& \sum_{R^1, R^2, R^3} \cW_{R^1R^2}\cW_{R^2R^3}\cW_{R^3R^1}
    (-1)^{|R^1|+|R^2|+|R^3|}q^{\half\sum_{i=1}^3k_{R^i}}\\
  && \times \exp\left\{ t_H(|R^1|+|R^2|+|R^3|) \right\}
\eea%
We need to study the following term
\bean%
 \cW_{R^1R^2}\cW_{R^2R^3}\cW_{R^3R^1}
    (-1)^d q^{\half\sum_{i=1}^3k_{R^i}} \label{eqn:degree of P2}
\eean%
under the condition $|R^1|+|R^2|+|R^3|=d$. Now
\bea%
 \cW_{\lambda\mu}=q^{-\frac{|\lambda|+|\mu|}{2}}
 S_\lambda(q^\alpha)S_\mu(q^{\lambda+\alpha})
\eea%
By (\cite{M} page 44 ex.1)
\bea%
 S_\lambda(q^\alpha)&=&q^{-n(\lambda)}\prod_{x\in\lambda}\frac{1}{1-q^{-h(x)}}
\eea%
We also have
\bea%
 S_\mu(q^{\lambda+\alpha})=\sum_{|\rho|=|\mu|}\frac{\chi_\mu(\rho)}{z_\rho}p_\rho
\eea%
Here we can see degree of $q$ in $S_\mu(q^{\lambda+\alpha})$ is
not greater than $|\lambda|\cdot |\mu|$ by writing down $p_\rho$
explicitly and equality holds when $\lambda=(k)$ for
$k=|\lambda|$. Hence the degree of $q$ in (\ref{eqn:degree of P2})
is no more than
$$
 |R^1||R^2|+|R^2||R^3|+|R^3||R^1|-d+\half\sum_{i=1}^3 k_{R^3}
$$
while
\bea%
 \half\sum_{i=1}^3 k_{R^i}
  &=& \sum_{i=1}^3 \left( n((R^i)^t)-n(R^i)
    \right)\\
  &\leq& \sum_{i=1}^3 n((R^i)^t)\\
  &\leq& \sum_{i=1}^3 \frac{|R^i|\cdot (|R^i|-1)}{2}
\eea%
The last inequality is obtained from $n(\lambda^t)\leq n((m))$
where $m=|\lambda|$. This can be proved in the following way. Let
$l=l(\lambda)$. If one subtracts one from $\lambda_l$ and adds one
to $\lambda_1$, one will find $n(\lambda^t)$ is increasing. So we
have the inequality. The equality holds when $\lambda=(m)$.

Let $|R^i|=x_i$. Then $\sum x_i=d$. In the above inequalities, all
equality hold when the corresponding partition is of length $1$.
So the $\deg P_d-1$ is
\bea%
 &&\sum_{i=1}^3 \frac{x_i\cdot (x_i-1)}{2} -d +
 (x_1x_2+x_2x_3+x_3x_1)\\
 &=& \frac{d^2-3d}{2}
\eea%
The coefficient of the highest degree in $P_\Sigma$ is the number
of the non-negative integer solutions of $x_1+x_2+x_3=d$
multiplied by $(-1)^d$, which is $\frac{(-1)^d}{2}(d+1)(d+2)$.

We do the similar calculation and get the degree of
$P_{(d,m)}^{(\bF_k)}$. The coefficient of the highest order term
is the number of non-negative integer solutions of
$x_1+x_3+kx_4=m,\ x_2+x_4=d$ multiplied by $(-1)^{dk}$, which is
$(-1)^{dk}\frac{(d+1)(m+1-kd)}{2}$.

Now we can easily verify the theorem by adjunction formula and
Hirzebruch-Riemann-Roch theorem.
\end{proof}

\appendix

\section{}

\begin{thm}\label{thm:number theory}
$\sum_{i=1}^n a_i = a$, $p$ is prime, $r\geq 1$. Then
\bea%
 \begin{pmatrix}
    p^r a \\
    p^r a_1, \cdots, p^r a_n
 \end{pmatrix}
 - \begin{pmatrix}
    p^{r-1}a \\
    p^{r-1} a_1, \cdots, p^{r-1} a_n
 \end{pmatrix}
 &\equiv& 0 \mod(p^{2r})
\eea%
\end{thm}

\begin{lemma}\label{lemma:p-r factor of n=2}
$p$ is prime, $r\geq 1$. Then
\bea%
 \begin{pmatrix}
    p^r a \\
    p^r b
 \end{pmatrix}
 - \begin{pmatrix}
    p^{r-1}a \\
    p^{r-1} b
 \end{pmatrix}
 &\equiv& 0 \mod(p^{2r})
\eea%
\end{lemma}

\begin{proof}
\bea%
 \frac{\begin{pmatrix}
    p^r a \\
    p^r b
 \end{pmatrix}}{\begin{pmatrix}
    p^{r-1}a \\
    p^{r-1} b
 \end{pmatrix}}
 &=&
 \frac{\prod_{k=1}^{p^rb}\frac{(a-b)p^r+k}{k}}{\prod_{k=1}^{p^{r-1}b}\frac{(a-b)p^{r-1}+k}{k}}\\
 &=& \prod_{\gcd(k,p)=1, 1\leq k \leq
    p^rb}\frac{(a-b)p^{r}+k}{k}\\
 &=& \prod_{\gcd(k,p)=1, 1\leq k \leq
    p^rb} \left( 1+\frac{a-b}{k}p^r \right) \\
 &\equiv& 1 + p^r(a-b)\sum_{\gcd(k,p)=1, 1\leq k \leq
    p^rb} \frac{1}{k} \mod(p^{2r})
\eea%
Let $A_p(n)=\sum_{1\leq k \leq n, \gcd(k,p)=1}\frac{1}{k}$.
Therefore the proof of the lemma is completed by showing
$A_p(p^rb)=\frac{p^rc}{d}$ for some $c,d$ such that $\gcd(d,p)=1$.

If $\gcd(k, p)=1$, there exist $\alpha_k$, $\beta_k$ such that
$$ \alpha_k k + \beta_k p^r =1 $$
Let $B_p(n)=\sum_{1\leq k \leq n, \gcd(k, p)=1}k$. Then by the
above formula
\bea%
 A_p(p^rb)
 &\equiv& b A_p(p^r) \mod(p^r)\\
 &\equiv& b B_p(p^r) \mod(p^r)
\eea%
And
\bea%
 B_p(p^r)
  &=& \sum_{k=1}^{p^r}k-p\sum_{k=1}^{p^{r-1}}k\\
  &=& \frac{p^r(p^r+1)}{2}-p\frac{p^{r-1}(p^{r-1}+1)}{2}\\
  &=& \frac{p^{2r-1}(p-1)}{2}
\eea%
$2r-1\geq r$ since $r\geq 1$. The proof is completed.
\end{proof}

Now, let's turn to the proof of Theorem \ref{thm:number theory}.
\begin{proof}
\bea%
  && \begin{pmatrix}
    p^r a \\
    p^r a_1, \cdots, p^r a_n
 \end{pmatrix}
 - \begin{pmatrix}
    p^{r-1}a \\
    p^{r-1} a_1, \cdots, p^{r-1} a_n
 \end{pmatrix}\\
 &=& \prod_{k=1}^n
    \begin{pmatrix}
        p^r(a-\sum_{i=1}^{k-1}a_i) \\
        p^r a_k
    \end{pmatrix}
  - \prod_{k=1}^n \begin{pmatrix}
        p^{r-1}(a-\sum_{i=1}^{k-1}a_i) \\
        p^{r-1} a_k
    \end{pmatrix}\\
 &\equiv& 0 \mod(p^r)
\eea%
In the last step, we have used Lemma \ref{lemma:p-r factor of
n=2}. The proof is then completed.
\end{proof}


\end{document}